\input amstex
\input amsppt.sty   
\input epsf.tex
\hsize 30pc
\vsize 47pc
\magnification=\magstep1
\def\nmb#1#2{#2}         
\def\cit#1#2{\ifx#1!\cite{#2}\else#2\fi} 
\def\totoc{}             
\def\idx{}               
\def\ign#1{}             
\redefine\o{\circ}

\define\al{\alpha}

\define\ep{\varepsilon}

\define\th{\theta}

\define\la{\lambda}
\define\rh{\rho}

\define\ph{\varphi}
\define\ch{\chi}

\define\De{\Delta}

\define\x{\times}
\define\p{\partial}
\def\today{\ifcase\month\or
 January\or February\or March\or April\or May\or June\or
 July\or August\or September\or October\or November\or December\fi
 \space\number\day, \number\year}
\topmatter
\title  Bounds on the multiplicity of eigenvalues for fixed membranes
\endtitle
\author  Thomas Hoffmann-Ostenhof \\
Peter W\. Michor  \\
Nikolai Nadirashvili
\endauthor
\affil
Erwin Schr\"odinger Institut f\"ur Mathematische Physik,
Boltzmanngasse 9, A-1090 Wien, Austria
\endaffil
\address
T\. Hoffmann-Ostenhof: Institut f\"ur Theoretische Chemie, 
Universit\"at Wien, W\"ahringer Stra{\ss}e 17, A-1090 Wien, Austria;
{\it and:}
Erwin Schr\"odinger Institut f\"ur Mathematische Physik,
Boltzmanngasse 9, A-1090 Wien, Austria
\endaddress
\email
thoffman\@esi.ac.at
\endemail
\address
P\. Michor:
Institut f\"ur Mathematik, Universit\"at Wien,
Strudlhofgasse 4, A-1090 Wien, Austria;
{\it and:}
Erwin Schr\"odinger Institut f\"ur Mathematische Physik,
Boltzmanngasse 9, A-1090 Wien, Austria
\endaddress
\email Peter.Michor\@esi.ac.at \endemail
\address 
N\. Nadirashvili:
Department of Mathematics, MIT, Cambridge, MA 02139, USA;
{\it currently:} The University of Chicago Department of Mathematics, 
5734 S. University Ave., 
Chicago, IL 60637, USA 
\endaddress
\email nikolai\@math.uchikago.edu \endemail
\dedicatory \enddedicatory
\thanks 
Work supported by the European
Union TMR grant FMRX-CT 96-0001 
\endthanks
\thanks
P.W.M\. was supported by `Fonds zur
F\"orderung der wissenschaftlichen                    
Forschung, Projekt P~10037~PHY'.
\endthanks
\abstract For a membrane in the plane the multiplicity of the $k$-th 
eigenvalue is known to be not greater than $2k-1$. 
Here we prove that it is 
actually not greater than $2k-3$, for $k\ge 3$.
\endabstract
\leftheadtext{\smc Hoffmann-Ostenhof, Michor, Nadirashvili}
\rightheadtext{\smc Bounds on the multiplicity of eigenvalues for membranes}

\endtopmatter

\document


\head\totoc\nmb0{1}. Introduction and Statement of the Result \endhead

Let $D\subset \Bbb R^2$ be a bounded domain with smooth boundary 
$\p D$. We consider the corresponding Dirichlet eigenvalue 
problem 
$$\cases
-\De u = \la_k u,\quad k=1,2,\dotsc, \quad  
\la_1<\la_2\le \la_3\le\la_4\le\dotso\\ 
u|\p D=0. 
\endcases\tag{\nmb.{1.1}}$$
For this problem we investigate the multiplicity of the eigenvalues 
$\la_j$, where $\la_j$ is said to have {\it multiplicity} 
$$
m(\la_k)=l\quad\text{ if }\quad 
\la_{k-1}<\la_k=\la_{k+1}=\dots=\la_j=\dots=\la_{k+l-1}<\la_{k+l}.
$$  
It is the dimension of the eigenspace 
$U(\la_k)=U(\la_{k+1})=\dots=U(\la_{k+l-1})$ of the eigenvalue 
$\la_k=\dots=\la_{k+l-1}$. 

Our goal is to find universal upper bounds for $m(\la_k)$.

From basic spectral theory it is known that $\la_1$ is simple. Cheng 
showed in a celebrated paper \cit!{4} that 
$$
m(\la_2)\le 3 
$$
for membranes and surfaces of genus 0. 
This is sharp for membranes, see \cit!{9}, where an example with 
$m(\la_2)=3$ is given; note that then also $m(\la_3)=3$.
There is very interesting work about $m(\la_2)$ for surfaces with 
genus $>0$, \cit!{2}, \cit!{5}, and \cit!{6}. 
It is known that in higher dimensions no 
universal bound to multiplicities can exist, \cit!{6}. About 10 years 
ago one of us \cit!{10} showed that 
$$
m(\la_k)\le 2k-1 
$$               
not only for the membrane case but also for Laplacians on surfaces 
with genus 0.

In a recent paper \cit!{8} it was shown for eigenvalues of 
Laplace Beltrami operators on smooth compact surfaces without 
boundary with genus 0 that 
$$
m(\la_k)\le 2k-3\quad\text{ for }\quad k\ge 3. 
$$
Here we prove the same result for the membrane case.

\proclaim{Theorem A}
Let $k\ge 3$. Then the multiplicity of the $k$-th eigenvalue $\la_k$ 
for the Dirichlet problem on $D$ satisfies 
$$
m(\la_k)\le 2k-3.
$$
\endproclaim

This will follow from the sharper Theorem B 
below.

\subhead\nmb.{1.2}. Remarks \endsubhead
The proof in \cit!{8} and the present one are quite different. 
Now we have a boundary and this requires a different approach though 
both proofs are based on a combination of Courant's nodal theorem, a 
suitable version of Euler's polyhedral theorem, and a detailed 
investigation of the zero sets of solutions $u$ of (\nmb!{1.1}). Here 
we have to investigate the zero sets near the boundary. 

The Laplacian in (\nmb!{1.1}) can be replaced by a strictly elliptic 
operator of second order in divergence form with smooth coefficients 
and one can also allow for a potential, so that theorem A 
holds also for the multiplicities of the eigenvalues of the following 
problem:
$$
\cases \Bigl(-\sum_{i,j=1}^2 
\frac{\p}{\p x_i}a_{i,j}\frac{\p}{\p x_j}+ 
V(x)\Bigr)u = \la_k u\quad & \text{ in }D\\
u=0 & \text{ in }\p D \endcases 
$$ 
We consider the principal symbol as the inverse of a Riemann 
metric on $D$ and use it to express 
angles etc\. in the proofs below.

Our result can be shown to carry over to the 
free membrane case, i.e\. $-\De u= \la_k u$ in $D$ with Neumann 
boundary conditions, but we do not go into 
details here. 

Probably one can relax the smoothness conditions considerably. It 
would be interesting to allow for unbounded regions, in particular 
for Schr\"odinger  equations in $\Bbb R^2$. 

For the membrane case there is an extensive literature on the 
asymptotics of eigenvalues. It is interesting to investigate 
the asymptotics of the following quantity:
$$
\Cal M(k) = \max \{m(\la_k): \text{ all membranes }D\}
$$

In section \nmb!{2} we shall recall some well known properties of 
eigenfunctions of membrane problems, state theorem B (a 
generalization of theorem A), give a suitable version 
of Euler's theorem on polyhedra, and prove that $m(\la_k)\le 2k-2$ 
for $k\ge 3$.
In section \nmb!{3} we complete the proof of theorem B. 

\head\totoc\nmb0{2}. Basics and the proof that $m(\la_k)\le 2k-2$ for 
$k\ge 3$. \endhead

\subhead\nmb.{2.1}. Nodal sets \endsubhead
Let $D$ be a bounded domain in $\Bbb R^2$ with smooth boundary which 
decomposes into connected components as 
$$
\p D = \bigcup_{i=1}^N (\p D)^i.
$$
We consider a solution $u$ of (\nmb!{1.1}), and define its 
\idx{\it nodal set} by
$$
\Cal N(u) := \overline{\{x\in D: u(x)=0\}}.
$$
It is well known (and follows from \nmb!{2.2}, \nmb!{2.4}, 
\nmb!{2.5}, and \nmb!{2.6} below) that:
\roster
\item $\Cal N(u)$ is a union of smoothly immersed circles in $D$ and 
       immersed arcs connecting points of $\p D$. Each of these 
       is called a \idx{\it nodal line}. Note that 
       self intersections are allowed. Maximal embedded pieces of 
       nodal lines will be called \idx{\it nodal arcs}. 
\item If $u(x)=0$ but $du(x)\ne0$ then $x$ lies on exactly one nodal 
       line and is no point of self intersection  of this nodal line.
\item If $u(x)=0, du(x)=0, \dots d^lu(x)=0$, but $d^{l+1}u(x)\ne0$ 
       then exactly $l+1$ nodal lines go through $x$ whose tangents 
       at $x$ dissect the full circle into $2l+2$ equal angles. In 
       particular the intersections are transversal. See \nmb!{2.4} 
       below.
\item If $x$ is a zero of order $l$ as in \therosteritem3 and lies in 
       $\p D$, then one of the nodal lines lies in 
       $\p D$. Here we use the fact the locally near $x$ the 
       eigenfunction $u$ may 
       be extended to the outside of $D$ and is there a solution of 
       the extended Laplace operator. 
\item Each component of the boundary is hit by an even number of 
       nodal lines (since $u$ changes sign at the nodal lines). 
\endroster

\subhead \nmb.{2.2} \endsubhead
The \idx{\it nodal domains} of $u$ are the connected components of 
$D\setminus \Cal N(u)$. We denote the number of nodal domains by 
$\mu(u)=\mu(\Cal N(u))$.  

\proclaim{Courant's Nodal Theorem} \cit!{3}
For each function $u$ in the eigenspace $U(\la_k)$ we have $\mu(u)\le k$.
\endproclaim

\subhead Remarks \endsubhead
It has been observed by Pleijel \cit!{11} that for the membrane case 
we have $\limsup \frac{\mu(u_k)}{k} < 1$ for any choice of 
$u_k\in U(\la_k)$. 
Or, in other words, Courant's nodal theorem is sharp only for 
finitely many eigenvalues. In view of this the following result is a 
substantial improvement of theorem A. 

\proclaim{Theorem B}
Let $U$ be a linear subspace of an eigenspace $U(\la)$
and let $1<l\in \Bbb N$ be such that for each 
$u\in U$ we have $\mu(u)\le l$. 
Then $\dim U \le \max\{3, 2l-3\}$.
\endproclaim

Note that now the inequality $m(\la_k)\le 2k-3$ for $k\ge 3$ in 
theorem A follows immediately from 
theorem B using Courant's nodal theorem 
\nmb!{2.2}. 

\subhead\nmb.{2.3}. Remarks \endsubhead
It would be interesting to investigate for which $l$ theorem B is 
sharp. 
For large $l$ it is not at all clear whether $\dim U$ can be 
estimated by something better than $O(l)$, such as $O(l^{1/2})$. 
There is a severe lack of examples with high multiplicities.

\proclaim{\nmb.{2.4}. Proposition} \cit!{1}, \cit!{4}
For an eigenfunction $u$ and $x_0\in D$ there exists an integer 
$n\ge0$ 
such that 
$$
u(x) = P_n(x-x_0) + O(|x-x_0|^{n+1})
$$
for a harmonic homogeneous polynomial $P_n\not\equiv 0$ of degree $n$.
\endproclaim

Actually for the membrane case we even have 
$$
u(x) = P_n(x-x_0) + P_{n+1}(x-x_0) + O(|x-x_0|^{n+2})
$$
for harmonic homogeneous polynomials $P_n\ne0$ and $P_{n+1}$ of degrees 
$n$ and $n+1$, respectively, but we shall not need this sharper 
result. 

Harmonic homogeneous polynomials $P_n$ of degree $n$ have a 
particularly simple representation in polar coordinates $(r,\th)$:
$$
P_n(r\cos\th,r\sin\th) = a r^n\cos(n\th) + br^n\sin(n\th).
$$
Obviously the set of zeros of such a $P_n$ consists of $n$ straight 
lines which meet at equal angles. 

\proclaim{\nmb.{2.5}. Proposition}
For an eigenfunction $u$ and $x_0\in\p D$ there exists a 
harmonic homogeneous polynomial $P_n\not\equiv 0$ of degree $n\ge1$
with
$$
u(x) = P_n(x-x_0) + O(|x-x_0|^{n+1})
$$
such that one of the nodal lines of $P_n$ is tangent to 
$\p D$ at $x_0$. 
\endproclaim

\demo{Proof}
This follows from the smoothness of $\p D$, see e.g.\ \cit!{7}.
\qed\enddemo

\subhead\nmb.{2.6} \endsubhead
Note that for $x_0\in D$ the leading 
harmonic homogeneous polynomial $P_n$ of an eigenfunction $u$ lies in 
the 2-dimensional (if $n>0$) vector space of all 
harmonic homogeneous polynomials of degree $n$, whereas for 
$x_0\in \p D$ it lies in the 1-dimensional subspace of those 
polynomials which vanish on the tangent $T_{x_0}(\p D)$.

\subhead\nmb.{2.7}. Definition \endsubhead
An \idx{\it (abstract) nodal set} is a set $\Cal N$ satisfying 
\nmb!{2.1}, \therosteritem1-\therosteritem5, 
where we do not require 
that it is the nodal set of an eigenfunction.

An \idx{\it isotopy} of nodal sets is 
a curve of nodal sets such that each immersed circle or arc moves 
along a smooth isotopy which respects nodal arcs. So intersection 
points can move but not change the multiplicity.

A \idx{\it nodal pattern} is an isotopy class of nodal sets. We shall 
often draw a clearly recognizable representative of a nodal pattern, 
see below.

\proclaim{\nmb.{2.8}. Proposition}
Let $\Cal N$ be an abstract nodal set in a domain $D$. 
Then we have
$$\align
\mu(\Cal N) = b_0 &(\Cal N\cup \p D) - b_0(\p D) 
     + \sum_{x\in \Cal N\cap D} (\nu(x)-1) + 
     \sum_{y\in\p D\cap \Cal N} \frac{\rh(y)}2 + 1,
     \quad\text{ where}\\
\mu(\Cal N) &= \text{ number of (nodal) domains of }D\setminus\Cal N,\\
b_0(\Cal N\cup \p D) &= \text{ number of connected components of 
}\Cal N\cup\p D,\\
b_0(\p D) &= \text{ number of connected components of 
the boundary }\p D,\\
\nu(x) &= \text{ number of nodal lines containing }x\in D,\\
\rh(y) &= \text{ number of nodal lines hitting }\p D\text{ in }y.
\endalign$$
Moreover,
$$
b_0(\Cal N\cup \p D) - b_0(\p D) 
     + \sum_{y\in\p D\cap \Cal N} \frac{\rh(y)}2 \ge 1,
$$
so that for $\Cal N\ne \emptyset$ we get 
$$
\mu(\Cal N) \ge
      \sum_{x\in \Cal N\cap D} (\nu(x)-1) + 2.
$$
\endproclaim

\demo{Proof}
Suppose that $\p D$ has $k$ components. Then for the Euler 
characteristic of $\overline D$ we have 
$\ch(\overline D)=2-b_0(\p D)$, which 
can be seen from a simple cell decomposition of $\bar D$. See e.g\. 
\cit!{12}.

We consider $\Cal N\subset \overline D$ and extend it to a cell 
decomposition of $\overline D$ with 
\roster
\item "$c_0$" many 0-cells, namely the points 
       $x\in \Cal N\subset D$ through which $\nu(x)>1$ nodal lines 
       pass (i.e.\ $2\nu(x)>2$ 1-cells emanate), the points 
       $y \in \Cal N\cap \p D$ in which  
       $\rh(y)>0$ nodal lines hit the boundary, and one extra point $z$ 
       on a smooth part of each of the 
       $b_{\Cal N}=b_0(\Cal N\cup \p D)$ connected components of 
       $\Cal N\cup \p D$.
\item "$c_1$" many 1-cells. First the 
       $\frac12(\sum_x2\nu(x)+\sum_y\rh(y))$ nodal arcs of $\Cal N$ 
       connecting the intersection points and boundary hitting points of 
       $\Cal N$. Second, the $\sum_y1$ smooth pieces of $\p D$ 
       lying between the hitting points $y$. 
       Moreover, $b_{\Cal N}$ more 1-cells, namely for each $z$ 
       either a smooth arc coming  
       from subdividing the smooth arc by choosing $z$, or a 1-cell 
       corresponding to a component of $\p D$ which is hit by 
       no nodal line. Finally, $b_{\Cal N}-1$ extra 1-cells connecting 
       the $b_{\Cal N}$ points $z$ on the components of 
       $\Cal N\cup \p D$ in a suitable way to each other or to 
       some of the $y$'s.
\item "$c_2$" many 2-cells. Note that none of the extra 1-cells 
       dissects a nodal domain, so we have $c_2=\mu(\Cal N)$.
\endroster
Thus for the Euler characteristic we have 
$$\multline
2-b_0(\p D) = \ch(\overline D) = c_0 - c_1 + c_2 \\
= \sum_x1 + \sum_y1 + b_{\Cal N} 
     - \sum_x\nu(x) -\sum_y\frac{\rh(y)}2 
     -\sum_y1 -b_{\Cal N} - (b_{\Cal N}-1)
     + \mu, 
\endmultline$$
and thus
$$\align
\mu &= \sum_x(\nu(x)-1)+\sum_y\frac{\rh(y)}2 + 
     b_0(\Cal N\cup \p D) - b_0(\p D) + 1.
\endalign$$
The last assertion follows by treating each connected component 
$(\p D)^i$ of the boundary separately, if $\Cal N\ne \emptyset$:
$$\multline
b_0(\Cal N\cup \p D) - b_0(\p D) 
     + \sum_{y\in\p D\cap \Cal N} \frac{\rh(y)}2 =\\
= b_0(\Cal N\cup \p D) 
     + \sum_i\biggl(\sum_{y\in(\p D)^i\cap \Cal N} 
     \frac{\rh(y)}2 -1\biggr) \ge 1. \qed
\endmultline$$
\enddemo

\proclaim{\nmb.{2.9}. Lemma}
Let $U$ be a linear subspace of dimension $m\ge1$ of an eigenspace 
$U(\la)$.

Then for each $x_0\in D$ there exists an eigenfunction 
$0\ne u\in U$ such that $d^lu(x_0)=0$ for  $0\le l<[m/2]$, 
where $[m/2]$ is the largest integer $\le m/2$.
If $d^{[m/2]}u(x_0)\ne 0$ we have
$$
u(x) = P_{[m/2]}(x-x_0) + O(|x-x_0|^{[m/2]+1}),
$$

On the boundary, for any choice of points 
$y_1,\dots,y_{m-1}\in \p D$ there exists 
an eigenfunction $0\ne u\in U$ such that at each $y_i$ 
at least one nodal line 
of $u$ hits $\p D$. Some points $y_i$ might coincide, in which 
case the corresponding number of nodal lines hit there.
\endproclaim

\demo{Proof}
This is linear algebra using \nmb!{2.4} and \nmb!{2.5}. 
\qed\enddemo

\subhead\nmb.{2.10} \endsubhead
Let $U(\la)$ be an $m$-dimensional eigenspace for an eigenvalue 
$\la$. Consider the unit sphere $S^{m-1}\subset U(\la)$ with 
respect to the $L^2$-inner product, say. 
For each $u\in S^{m-1}$ we may consider its 
nodal set $\Cal N(u)$. We get a disjoint decomposition of $S^{m-1}$ 
(actually of $P^{m-1}(\Bbb R)$) according to the nodal patters. This 
should actually be a stratification into smooth manifolds. 

\proclaim{\nmb.{2.11}. Lemma}
Let $\ph:\Bbb R\to S^{m-1}\subset U(\la)$ be smooth.

Then for each multiindex $\al$ we have 
$$
\sup_{y\in\overline D} |(\p_x)^\al(\ph_{t}-\ph_s)| \le C_\al|t-s|
$$
for some constant $C_\al$.
\endproclaim

\demo{Proof}
This follows from the assumptions that all data are smooth.
\qed\enddemo

\subhead\nmb.{2.12} \endsubhead
Already from lemma \nmb!{2.9} we can prove the main result from 
\cit!{10} that $m(\la_k)\le 2k-1$ as 
follows: Suppose that $m(\la_k)\ge 2k$ and pick $x_0\in D$, then 
there is an eigenfunction $u\in U(\la_k)$ with $\nu_u(x_0)=k$, by 
lemma \nmb!{2.9}. Hence by lemma \nmb!{2.8} we get
$\mu(u) \ge 2-1 + k-1 + 1 = k+1$,
a contradiction to Courant's nodal theorem \nmb!{2.2}. 

Actually we even proved: 
{\it If $U$ is a linear subspace of an eigenspace 
$U(\la)$ and if $\sup\{\mu(u):u\in U\}= l>1$, then 
$\dim(U)\le 2l-1$.}

But we can do even better:

\proclaim{\nmb.{2.13}. Lemma}
Let $U$ be a linear subspace of an eigenspace $U(\la)$ and suppose 
that $\sup\{\mu(u):u\in U\setminus 0\}= l>1$.

Then we have $\dim(U)\le 2l-2$.
\endproclaim

\demo{Proof}
Assume for contradiction that there is some $U\subseteq U(\la)$ with 
$\dim(U)=2l-1$. We first assume that $D$ is simply connected. We 
pick $y,z\in\p D$. By lemma \nmb!{2.9} there exists an 
eigenfunction $u=u_{y,z}\in U$ such that 
$\rh_u(y)=2l-3$ and $\rh_u(z)=1$. By lemma \nmb!{2.8} and by 
\nmb!{2.2} we get 
$\mu(u)=l$, $\nu(x)-1=0$ for all $x\in \Cal N\cap D$, and 
$b_0(\Cal N\cup \p D)=1$. 

Consider now $\Cal N(u)$: there are $2l-3$ nodal lines emanating from 
$y$ and one from $z$, and there is no intersection point in $D$.
Let $\tilde {\Cal N}=\Cal N(u)\setminus \p D$, which consists 
of one smooth arc with endpoints $y$ and $z$, and of $l-2$ 
non-intersecting loops starting at $y$. For $l=4$ e.g., we have one 
of the following 5 nodal patterns:
\vskip 3mm 
\centerline{\epsfxsize=8cm\epsfbox{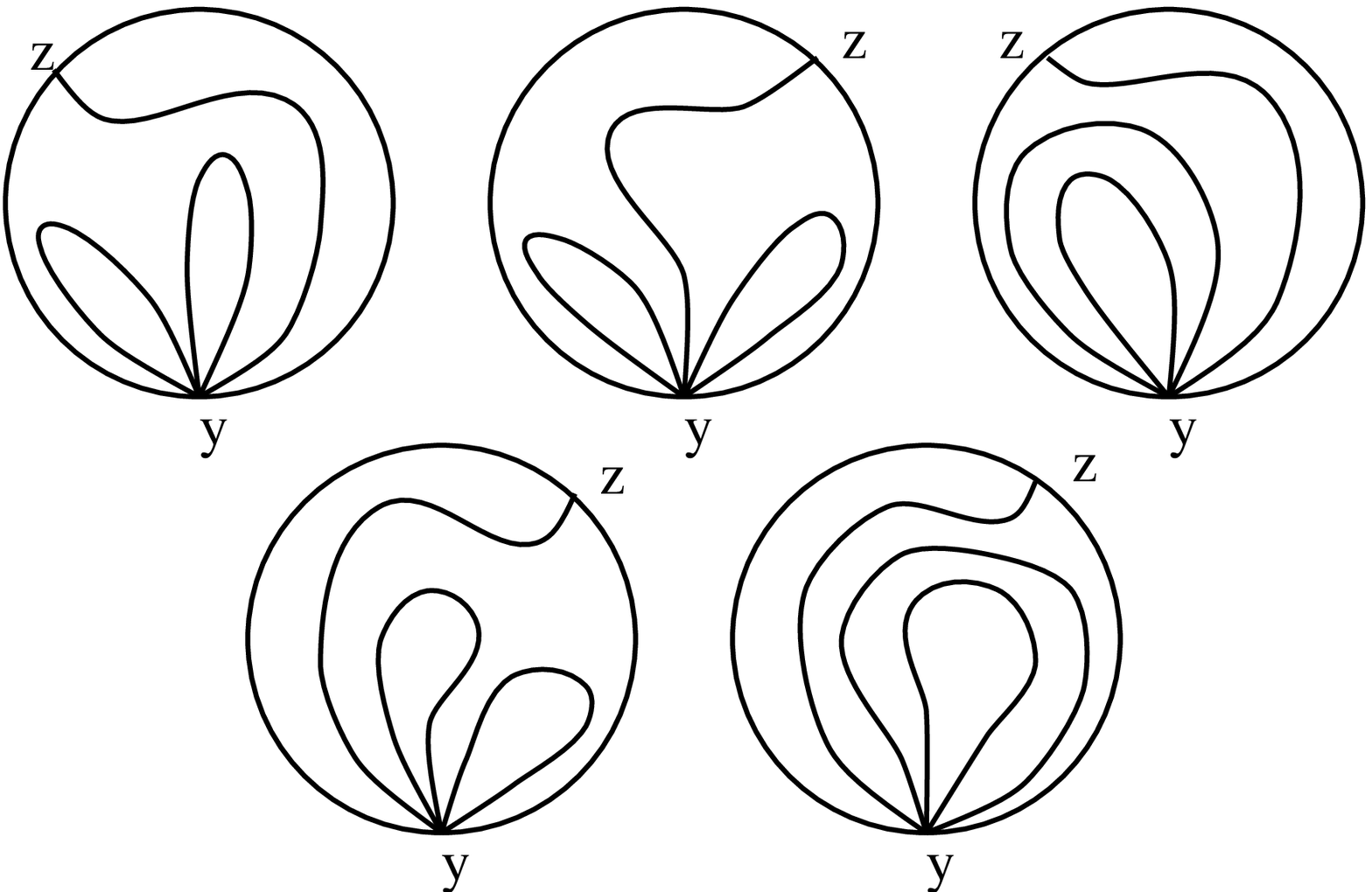}}
\medskip

We also note that for given $y,z\in \p D$ the eigenfunction 
$u_{y,z}$ is unique up to multiplication by a constant. 
Indeed, if there are two linearly 
independent eigenfunctions $u^1_{y,z}$ and $u^2_{y,z}$, 
then by \nmb!{2.5} there is an eigenfunction 
$v\in \text{span}(u^1,u^2)$ with $\rh_v(z)\ge 2$. Via \nmb!{2.8} (see 
\nmb!{2.1}.\therosteritem5)
we get as above for $l>1$ a contradiction to 
$\sup\{\mu(u):u\in U\}= l$.

Now we move $z$ towards $y$, once clockwise and once 
anticlockwise. Since we work at the maximal number of nodal domains, 
$\mu(u_{y,z})=l$, no additional intersection 
points in $D$ nor additional hitting points in $\p D$ may 
appear during these moves. Hence the arc from $y$ to $z$ will 
eventually become a loop as $z\to y$. But the limit nodal 
patterns differ, which is obvious from the figures above. For example 
\vskip 3mm 
\centerline{\epsfxsize=2.2cm\epsfbox{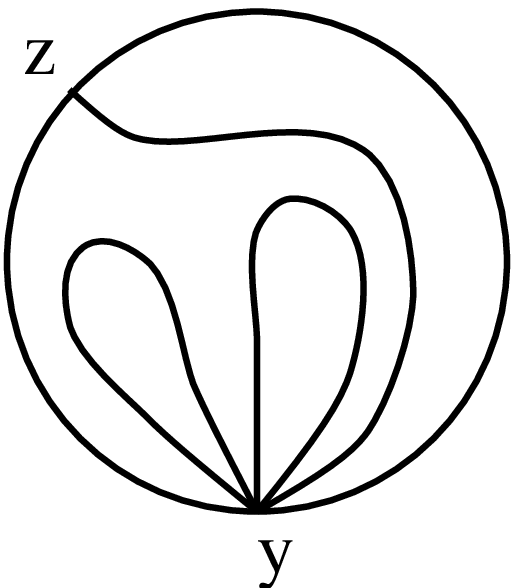}}
\medskip
will eventually tend to: 
\vskip 3mm 
\centerline{\epsfxsize=5cm\epsfbox{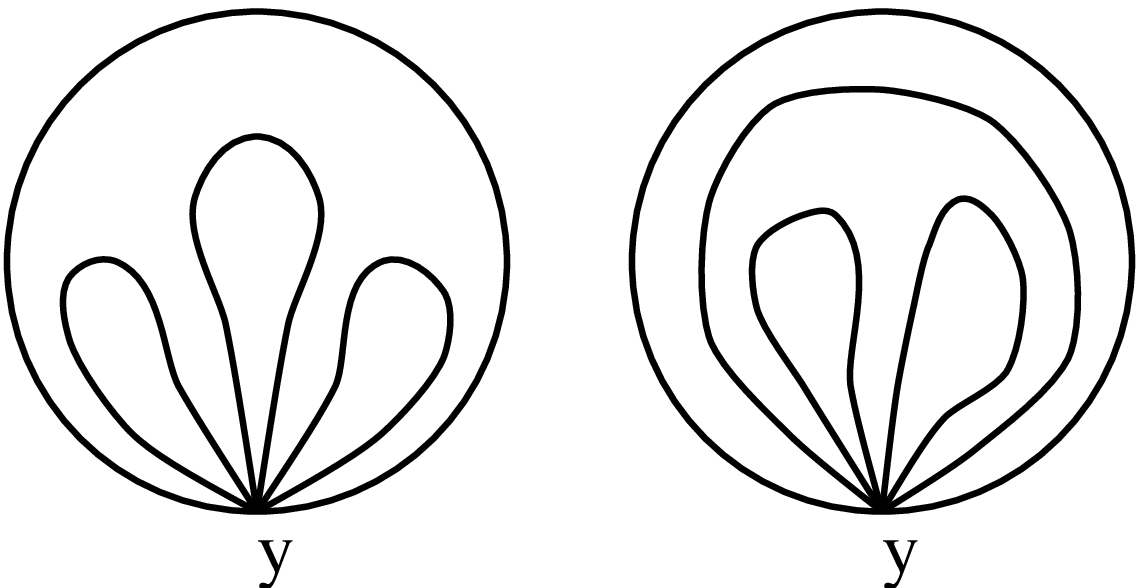}}
\hfill(clockwise) \hfill (anticlockwise) \hfill
\break
Hence there are two linearly independent functions with $2l-2$ nodal 
lines hitting at $y$. So again there is a function $v$ in their 
span with $2l-1$ nodal lines hitting at $y$, by \nmb!{2.5} a 
contradiction to $\mu(v)\le l$. 

The case of non simply connected domains is similar. We pick $y$ 
and $z$ as above on the outer component of the boundary and we 
proceed as above. Any other component of the boundary can be hit at 
most twice by nodal lines, according to \nmb!{2.8}, 
since we work at the maximal number $l$ of 
nodal domains. 
\qed\enddemo

\head\totoc\nmb0{3}. Proof of theorem B
\endhead

\subhead\nmb.{3.1}
\endsubhead
Let $U$ be a linear subspace of an eigenspace $U(\la)$ and let
$\sup\{\mu(u):u\in U\}=l\ge3$. 
We already know from \nmb!{2.13} that then $\dim(U)\le 2l-2$, so let 
us assume for contradiction that $\dim (U)=2l-2$, throughout this 
section. In all our constructions below we will use only 
eigenfunctions $u$ for which the number of nodal domains has to be 
maximal, i.e\. $\mu(u)=l$.

Before going into details we sketch the main ideas of the proof. We 
shall show that for each $x\in D$ there is a unique (up to 
multiplication by a constant) eigenfunction $u_x$ which vanishes of 
order $l-1$ at $x$. On the boundary, for each $y\in \p D$ there 
exists also a unique function $u_y$ (up to a multiplicative constant) which 
vanishes of order $\ge 2l-2$. We will show that these combine to 
a continuous mapping from $\overline D$ to the projective space $P(U)$. 
Then we shall use a winding number argument to get a contradiction. 

We start by giving a list of possible nodal 
patterns which hit $\p D$ in two points $x,y$ with 
$\rh(y)=2l-4$ or $2l-3$ and $\rh(x)=2$ or $\rh(x)=1$: We give all 
possible configurations at $x$, but just a sample of those possible 
at $y$, and we assume that $D$ is simply connected. All of these 
configurations look similar. We split each nodal pattern into 
two parts, namely into `the loops hitting the boundary only at $y$', 
and the rest, which can be either one nodal arc from $y$ to $x$, or
a loop hitting only at $x$, called a \idx{\it drop}, or 
two nodal arcs from $y$ to $x$, called a \idx{\it banana}.
\vskip 3mm 
\centerline{\epsfxsize=9.5cm\epsfbox{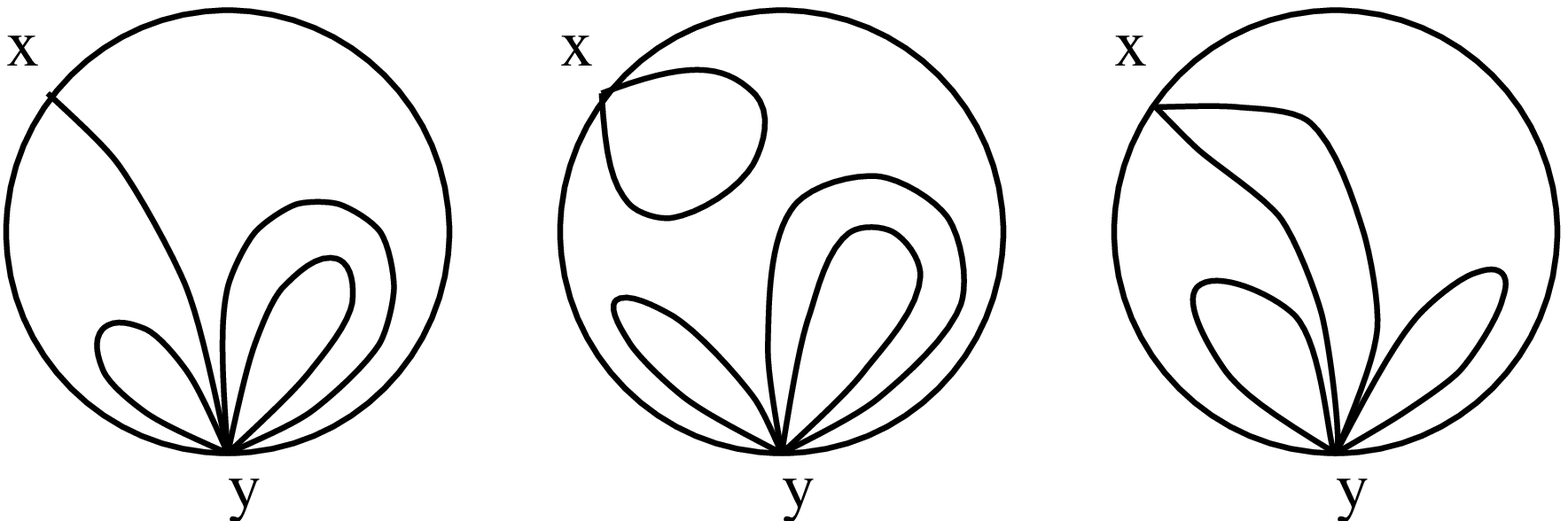}}
\smallpagebreak

All pictures and most arguments below will be given in the case that 
$D$ is simply connected. But since we always work with eigenfunctions 
which have the maximal number of nodal domains allowed by \nmb!{2.8}, 
everything remains valid in the non simply connected case: Then we 
have further boundary components each of which can be hit at most 
twice by nodal lines: otherwise we get too many 
nodal domains. Furthermore all boundary components are equivalent for 
our arguments (put $D$ into $S^2$), and we shall treat each of them 
separately.

\subhead\nmb.{3.2} \endsubhead
By \nmb!{2.9}, for each $y\in \p D$ there is a function 
$u_y\in U$ 
such that at least $2l-3$ nodal lines hit $\p D$ at $y$. 
The nodal pattern of $u_y$ consists thus of loops at $y$ and one or 
no nodal line from $y$ to another point $x\ne y$, in the simply 
connected case. In the general case it is similar with the changes 
described at the end of \nmb!{3.1}. 

\proclaim{Lemma} There are no
two linearly independent functions $u,v\in U$ such that $\rh_u(y)=2l-3$ 
and $\rh_v(y)=2l-2$. Moreover,
the set of points $y\in \p D$ where there 
exists a $u$ such that $\Cal N(u)$ has just $l-1$ loops at $y$, i.e\. 
$\rh_y(u)=2l-2$, is discrete. 
\endproclaim

\demo{Proof}
The nodal pattern $\Cal N(u)$ consists of $l-2$ loops at $y$ and one 
nodal line  
from $y$ to some point $x\ne y$, whereas $\Cal N(v)$ consists only of 
$l-1$ loops at $y$.
By a linear combination of $u$ and 
$v$ we may move $x$ anticlockwise or clockwise to $y$ 
and produce a function $w\in U$ such that $\Cal N(w)$ consists of 
loops at 
$y$ which is different from $\Cal N(v)$ (see \nmb!{2.13} for a 
similar argument). 
But by \nmb!{2.5} the leading terms of $w$ and $v$ at $y$ are multiples 
of the same harmonic polynomial, so a suitable linear combination of 
$w$ and $v$ has a zero of order at least $2l$ at $y$ which 
contradicts our assumption on $U$.

For the proof of the second assertion,
suppose that there is an open arc $I$ in $\p D$ such that for 
each $y\in I$ there exists $u_y\in U$ with $\rh_{u_y}(y)=2l-2$. 
Then by the argument above for the first assertion, $u_y$ is uniquely 
determined by $y$ and the mapping $I\ni y\mapsto u_y$ into the 
projective space of $U$ is 
smooth, since $u_y$ is given (up to a multiplicative constant) by 
solving a system of linear equations of maximal rank.
Let $y(t)=(y_1(t),y_2(t))$ be a unit speed parametrisation of $I$. 
Then near $y(t)$ the eigenfunction $u_{y(t)}$ with $2l-2$ nodal lines 
hitting $y(t)$ can be written as
$$
u_{y(t)}(x) = f(t)\Bigl(c_1(t) P^1_{2l-1}(x-y(t))
     + c_2(t) P^2_{2l-1}(x-y(t))\Bigr) + O(|x-y(t)|^{2l})
$$
where $P^1_{2l-1} = r^{2l-1}\cos((2l-1)\th)$ 
and $P^2_{2l-1} = r^{2l-1}\sin((2l-1)\th)$
span the 2-dimensional space of harmonic polynomials of degree 
$2l-1$,
where $f(t)$ is a normalizing function, and where 
$(c_1(t),c_2(t))\in S^1$  
is chosen in such a way that the leading term vanishes along 
the tangent $T_{y(t)}(\p D)$ spanned by $\dot y(t)$. 
We have $\p_t u_{y(t)}\in U$, and we compute this for a point 
$t_0$ where we may assume without loss that $\dot y(t_0)=(1,0)$, so that 
$c_1(t_0)=0$ and $c_2(t_0)=1$. Then 
$$
(\p_t u_{y(t)})|_{t=t_0}(x) 
     = f(t_0) \frac{\p P^2_{2l-1}}{\p x_1}(x-y(t_0)) 
     + O(|x-y(t_0)|^{2l-1}),
$$
where the leading term of order $2l-2$ does not vanish.
So in $\Cal N(\p_t(u_{y(t)}|_{t=t_0}))$ we have $2l-3$ nodal lines 
hitting $y(t_0)$, in 
contradiction to the first assertion of the lemma. 
\qed\enddemo

\proclaim{\nmb.{3.3}. Lemma} For each $y\in\p D$ there 
exists a unique (up to a nonzero constant) function $u_y\in U$ such 
that $\rh_{u_y}(y)\ge 2l-3$. Moreover, $y\mapsto u_y$ is a smooth map 
from $\p D$ into the projective space $P(U)$ of $U$.   
Put
$$\align
(\p D)_{2l-3} :&= \{y\in \p D: \rh_{u_y}(y)=2l-3\},\\
(\p D)_{2l-2} :&= \{y\in \p D: \rh_{u_y}(y)=2l-2\},\\
\endalign$$
then we have the disjoint union 
$\p D= (\p D)_{2l-3}\cup (\p D)_{2l-2}$, 
where $(\p D)_{2l-2}$ is discrete. 
Thus $(\p D)_{2l-3}$ is a union of open arcs and the nodal 
pattern of $u_y$ is constant for $y$ in one of these arcs. 
\endproclaim

Below is a 
list of such nodal patterns. 
Note that if $y$ moves to one of the endpoints $y_i$ of an interval of 
$(\p D)_{2l-3}$ then the last hitting point $z(y)$ of $u_y$ has to 
move towards $y_i$ too.
\vskip 3mm 
\centerline{\epsfxsize=9cm\epsfbox{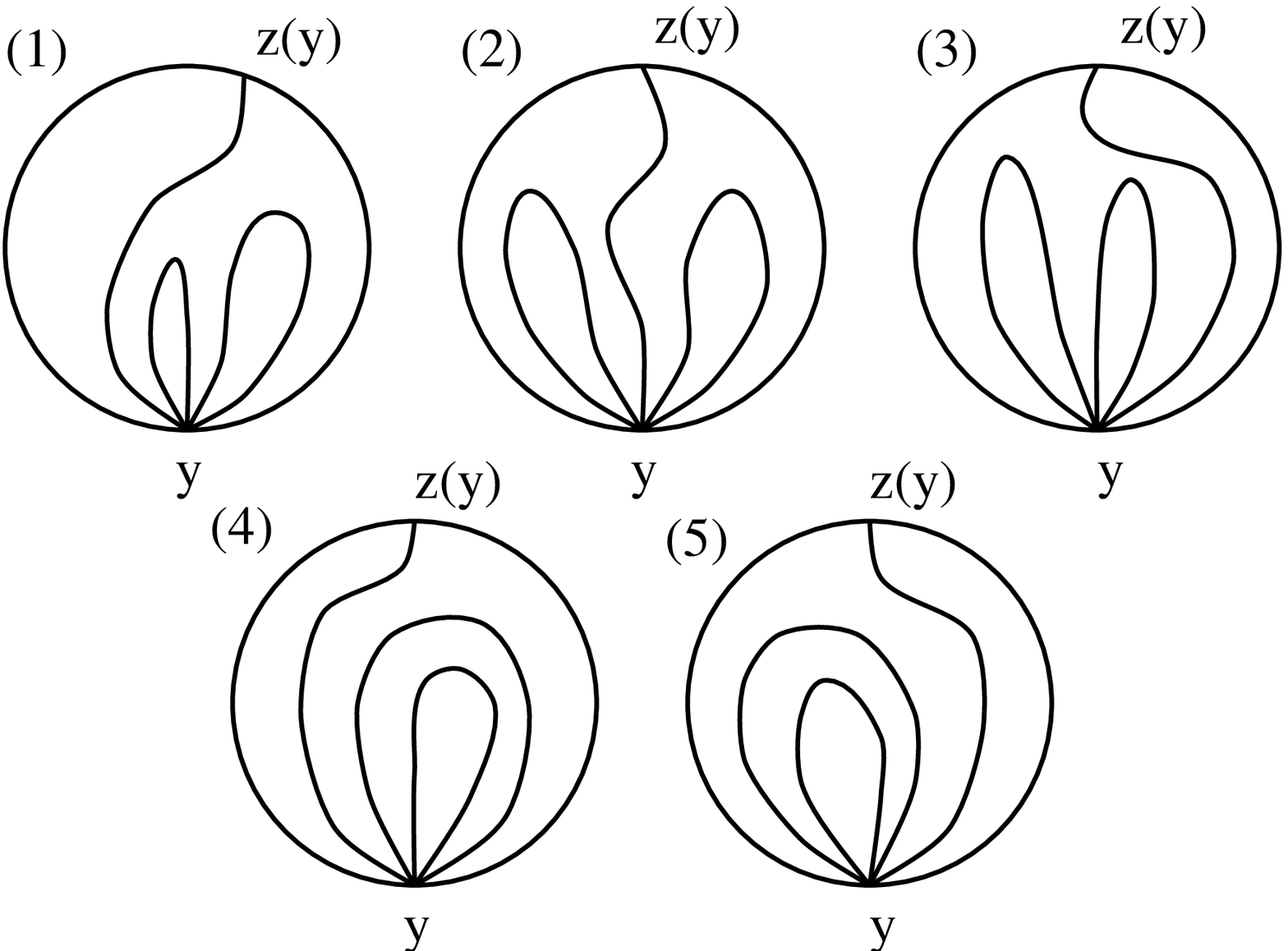}}

\demo{Proof}
If there are two linearly independent functions with $2l-3$ nodal 
lines hitting at $y$, a suitable linear combination has leading term 
of order one higher, so $2l-2$ nodal lines hitting at $y$, thus 
$y\in (\p D)_{2l-3}\cap (\p D)_{2l-2}$ which is empty by 
\nmb!{3.2}. If there are two linearly independent functions with 
$\rh(y)=2l-2$, a suitable linear combination has $\rh(y)=2l-1$ and 
thus too many nodal domains. 

The map $y\mapsto u_y$ is smooth by uniqueness and \nmb!{2.9}, since 
we solve there a linear system which has maximal rank by uniqueness.  

The rest is clear from \nmb!{3.2}.
\qed\enddemo

\proclaim{\nmb.{3.4}. Lemma}
For each $y\in (\p D)_{2l-3}$ and $s\in\p D$ there 
exist a function $v_{y,s}\in U$, unique up to a constant, with 
$\rh_{v_{y,s}}(y)\ge 2l-4$ and $\rh_{v_{y,s}}(s)\ge 1$. Moreover, 
$v_{y,z(y)}=u_y$, and $(y,s)\mapsto v_{y,s}$ is a smooth mapping from 
$(\p D)_{2l-3}\x \p D$ to $P(U)$.
\endproclaim

\demo{Proof}
Existence follows from \nmb!{2.9}. As explained in the proof of 
\nmb!{3.3}, smoothness follows from uniqueness, which we prove now. 

Suppose that $s\ne z(y)$ and that there are two linearly independent 
functions of this kind. Then a suitable linear combination has 
$\rh(y)=2l-3$ and $\rh(s)=1$ which contradicts the uniqueness of 
$u_y$. 

If $s=z(y)$ and there exists a second function $v_{y,z(y)}$ 
which is linearly independent of $u_y$, then 
the nodal pattern of 
$v_{y,z(y)}$ is typically of the form
\vskip 3mm 
\centerline{\epsfxsize=7cm\epsfbox{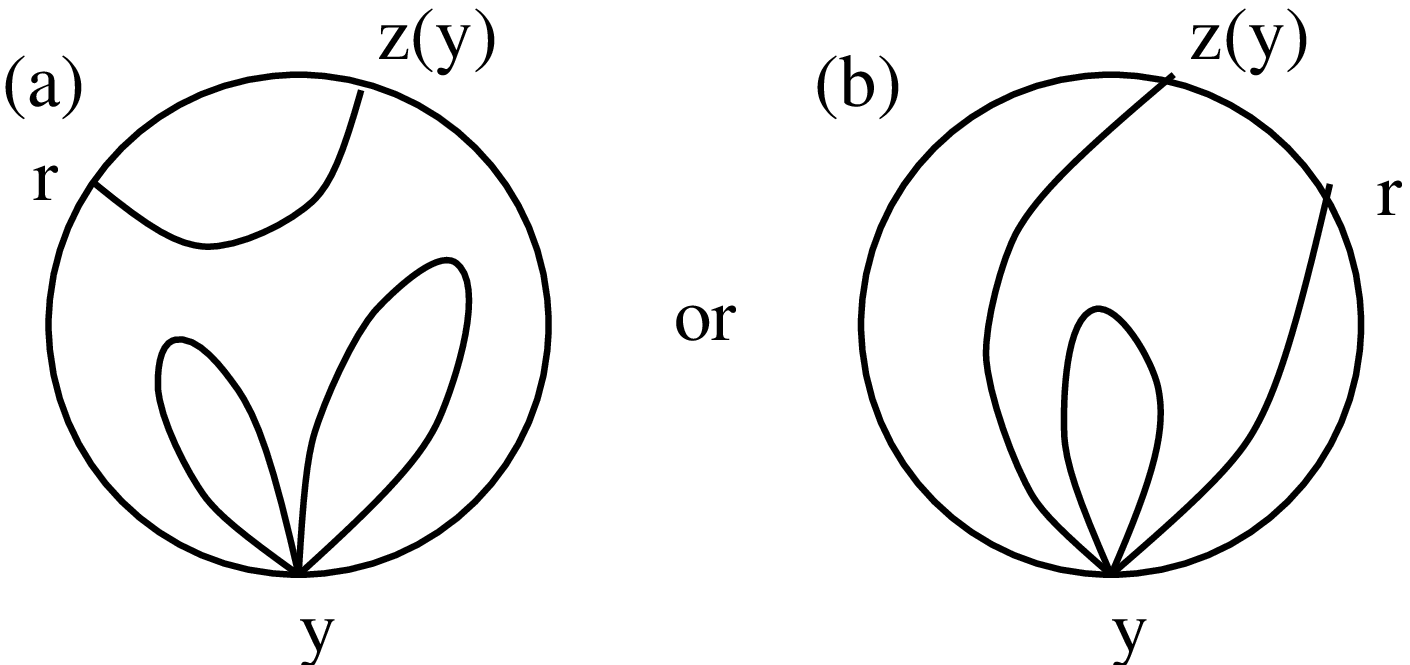}}
\medpagebreak
\noindent
since the third possibility (a nodal line from $z(y)$ to $y$ and all 
others loops at $y$) contradicts the uniqueness of $u_y$. 
In case (a) above, we arrange the signs of the functions 
$v_{y,z(y)}$ and $u_y$ as follows
\vskip 3mm 
\centerline{\epsfxsize=7cm\epsfbox{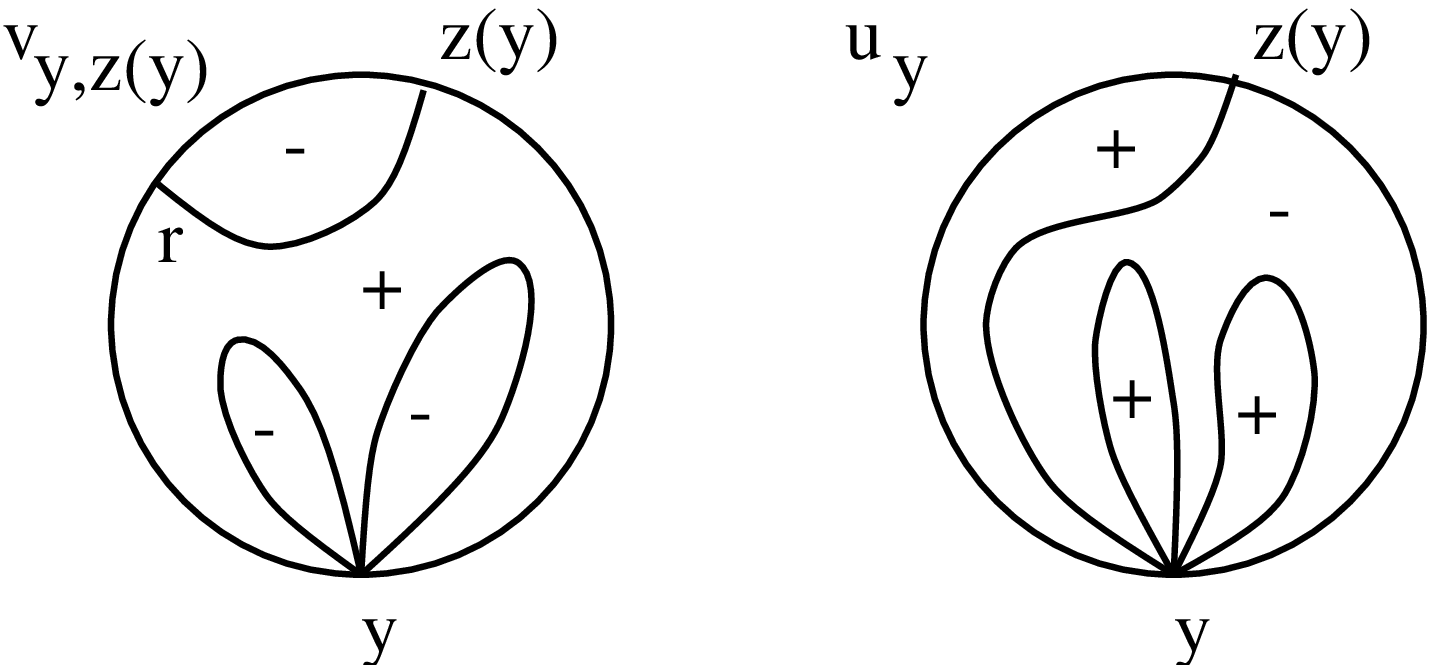}}
\medpagebreak
\noindent
and look at the nodal pattern of 
$w_t:=tu_y +(1-t)v_{y,z(y)}$. This can be viewed as follows: put the two 
drawings above each other and start at $t=0$, at 
$\Cal N(v_{y,z(y)})$. With growing $t$, domains where both functions 
are negative or where both functions are positive grow, whereas 
domains with mixed signs shrink. Thus the hitting point $r$ moves 
towards $z(y)$ and we get eventually, at some $0<t_1<1$, a nodal domain 
with a drop at $z(y)$. Further increasing $t$ this drop has to open 
but one nodal line has to stay at $z(y)$.
If it opens to the right we can never get the nodal pattern of $u_y$. 
If it opens to left the nodal line would eventually get to the point 
$r$ again, but then we would have two linearly independent functions $w_0$ 
and $w_{t_2}$ with $0<t_1<t_2<1$ with $2l-4$ nodal lines hitting at $y$ 
and one each at $r$ and $z(y)$. By a suitable linear combination we 
can then produce a function with $\rh(y)\ge 2l-3$, and nodal lines 
hitting at $r$ and $z(y)$, contradicting the uniqueness of $u_y$. 

Case (b) above is similar. These are the most complicated cases. 
Similar but more obvious methods apply 
if $r$ and $z(y)$ change position. If the nodal pattern of 
$u_y$ is different, such that the nodal line from $y$ to $z(y)$ 
has loops to the left and to the right, then
the argument is even easier.
\qed\enddemo

\proclaim{\nmb.{3.5}. Lemma}
Let the open arc $I$ be a connected component of 
$(\p D)_{2l-3}$ with endpoints $y_1$ and $y_2$ as in the 
drawing below. Let $y\in I$ and let $z(y)$ be the ultimate hitting 
point of $u_y$. Then $z(y)\notin I$. 
\endproclaim

\demo{Proof} Let us assume for contradiction that $z\in I$.
\vskip 3mm 
\centerline{\epsfxsize=3.0cm\epsfbox{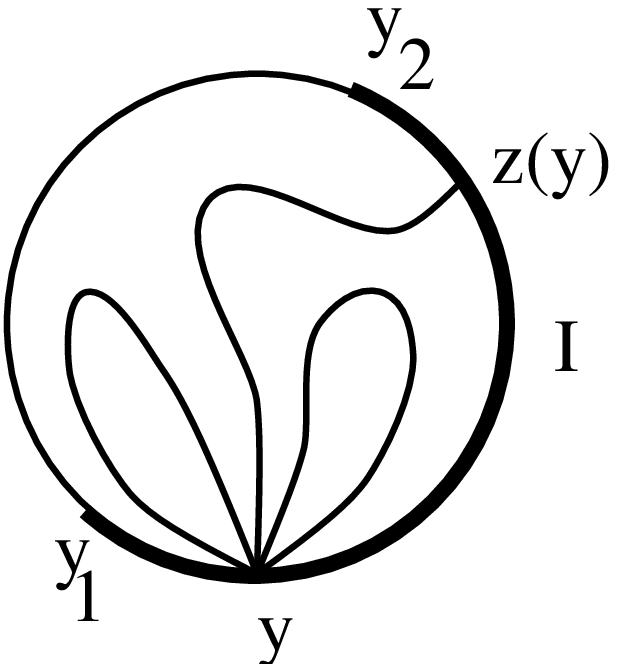}}
\medpagebreak
\noindent
We consider the function $v_{y,s}$ from lemma \nmb!{3.4}. Then 
$v_{y,z(y)}=u_y$, and we move $s$ inside $I$ towards $y$ (down right 
in the drawing). Then one of the nodal lines hitting at $y$ must move 
away from $y$ (otherwise we get a contradiction to the uniqueness of 
$u_y$): If this is the leftmost we must get eventually 
\vskip 3mm 
\centerline{\epsfxsize=3.0cm\epsfbox{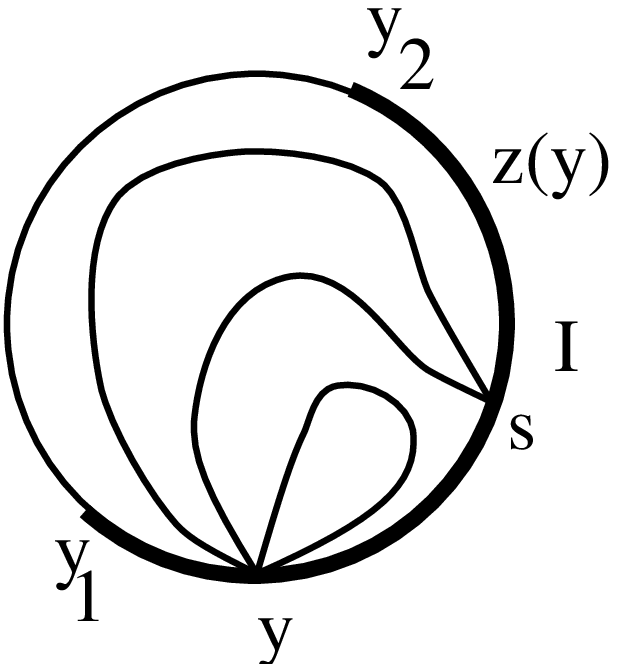}}
\medpagebreak
\noindent
because $s$ cannot move to $y$ for this would lead to a nodal type 
different from that of $u_y$. Call the corresponding function $v$. 
The functions $v$ and $u_y$ cannot coexist since in their span there 
is a function with nodal pattern
\vskip 3mm 
\centerline{\epsfxsize=3.0cm\epsfbox{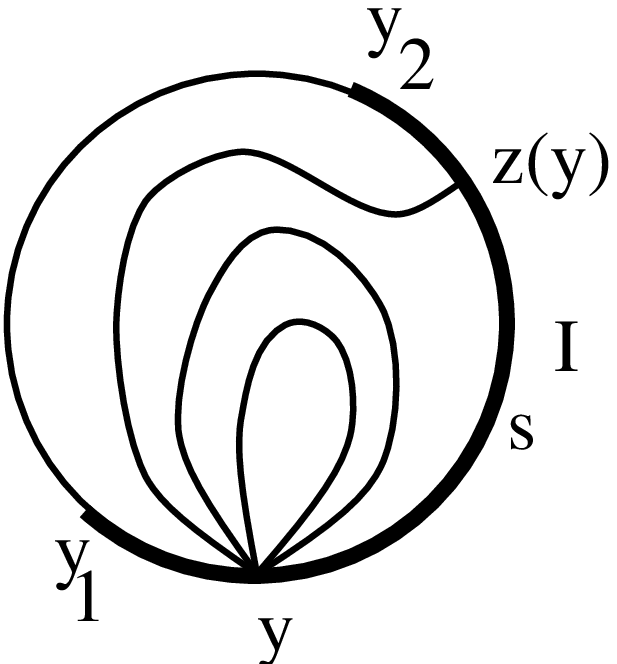}}
\medpagebreak
\noindent
which contradicts the nodal type valid on $I$. 

Hence the rightmost nodal arc must move away from $y$. If it moves 
down to $y$ we have already a contradiction. Thus it must 
eventually hit the downcoming $s$ at $s_1$ so that we 
have the the following nodal domain. 
\vskip 3mm 
\centerline{\epsfxsize=3.0cm\epsfbox{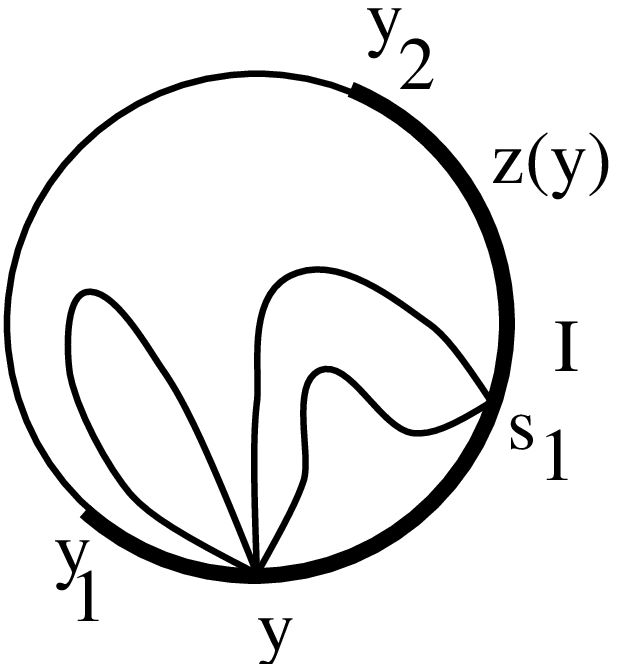}}
\medpagebreak
\noindent
Then we move $y$ towards $s_1$ and consider $\Cal N(v_{y,s_1})$. One
of the nodal lines hitting at $s_1$ must move away before $y$ hits 
$s_1$ since there is no point of $(\p D)_{2l-2}$ in between, 
and it must move eventually towards $y$ so that at some $x_1\in I$ between 
$y$ and $s_1$ we get $v_{x_1,s_1}=u_{x_1}$, since the nodal type of 
$u_y$ is constant in $I$. We have then the same situation as at the 
beginning, and we start to move again $s$ from $s_1$ to $x_1$ and consider 
$\Cal N(v_{x_1,s})$, and so on. We get a sequence of points $x_i$ and 
$s_i=z(x_i)$ in $I$ which move together. Even if they accumulate, at 
any accumulation point we must have the same nodal type, 
and we can continue the procedure. 
So we assume that finally $x_i\to x$ in $I$ and also $s_i\to x$.
But since $s_i=z(x_i)$ we finally get $x=z(x)$ so that 
$x\in (\p D)_{2l-2}$, a contradiction. 
\qed\enddemo

\proclaim{\nmb.{3.6}. Lemma}
Let the open arc $I$ be a connected component of 
$(\p D)_{2l-3}$ with endpoints $y_1$ and $y_2$ as in the 
drawing below. Let $y\in I$ and let $z(y)$ be the last hitting point of 
$u_y$. 

Then the following holds:
If $y$ moves clockwise to $y_1$ then $z(y)$ moves 
anticlockwise to $y_1$. If $y$ moves anticlockwise to $y_2$ then 
$z(y)$ moves clockwise to $y_2$. In particular the nodal patterns 
$\Cal N(u_{y_1})$ and $\Cal N(u_{y_2})$ are different. Moreover 
$(\p D)_{2l-3}$ consists of only finitely many open arcs, and 
$(\p D)_{2l-2}$ is a finite subset of $\p D$.
\endproclaim
Here is a sample of the nodal patterns of $u_x$ for $x=y_1$, 
$x\in I$, and $x=y_2$.
\vskip 3mm
\centerline{\epsfxsize=8cm\epsfbox{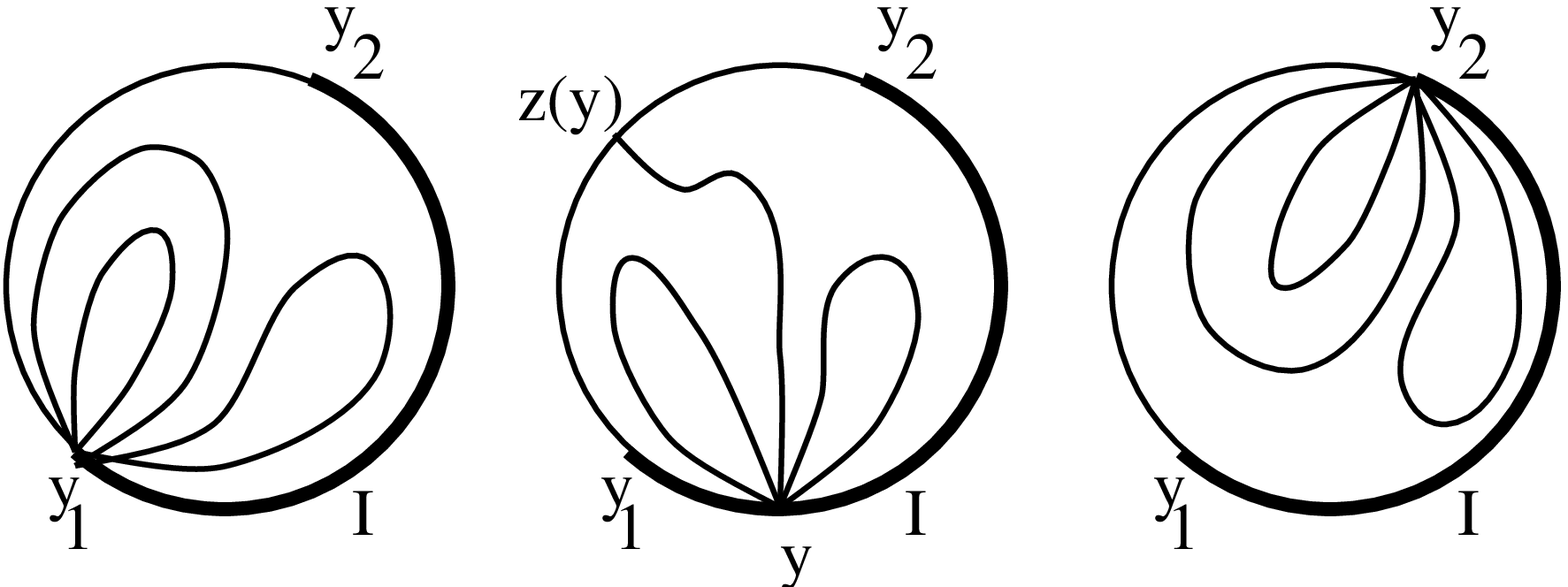}}

\demo{Proof}
Since $u_y$ depends smoothly on $y$, $z(y)$ also depends smoothly on 
$y$. If $y$ moves to $y_1$, $z(y)$ has to go to $y_1$ too since 
$\Cal N(u_{y_1})$ consists of loops at $y_1$ only. But it cannot come 
through $I$, by \nmb!{3.5}, so it must come from the outside.

Looking at the possibilities for $\Cal N(u_{y_i})$ one sees that 
there must be different nodal patterns at both ends of the arc $I$. 

If $(\p D)_{2l-2}$ were not finite, its points would accumulate 
at $y_0$, say. But then for $y_0$ there have to be two functions 
$u_{y_0}$ with different nodal patterns, a contradiction.
\qed\enddemo

\proclaim{\nmb.{3.7}. Lemma}
For each point $x\in D$ there exists a function $u_x\in U$, unique up 
to a multiplicative constant, such that $\nu_{u_x}(x)=l-1$. Moreover, 
$x\mapsto u_x$ induces a smooth mapping $D\to P(U)$ into the 
projective space of $U$.
\endproclaim

\demo{Proof}
Existence of $u_x$ follows from lemma \nmb!{2.9} and $\dim U=2l-2$.
From \nmb!{2.8} we see that 
$l\ge\mu(u_x)\ge \sum_{z\in D}(\nu_{u_x}(z)-1) +2\ge l$ so that $x$ 
is the only intersection point of $\Cal N(u_x)$ in $D$. 
{\it At most two nodal lines can connect $x$ to the (outer) boundary.}

If there are two linearly independent functions with the properties 
of $u_x$, we may choose 
functions $u_0$ and $u_1$ in their span
such that in local Riemannian polar coordinates 
$(r,\th)$ centered at $x_0$ we have
$$
u_0= r^{l-1}\cos((l-1)\th) + O(r^{l}),\quad
u_1= r^{l-1}\sin((l-1)\th) + O(r^{l}).
$$
Let $v_\al:=\cos(\al).u_1 +\sin(\al).u_2$, then $v_0=u_0$, and the 
regular $(l-1)$-gon consisting of the tangents to the 
$l-1$ nodal lines through $x_0$ in $T_{x_0}M$ 
rotates with $\al$ and is the same at the angle 
$\al=\pi/(l-1)$.
Thus $v_{\pi/(l-1)}$ has the same leading term at $x_0$ as 
$-u_0$, thus $v_{\pi/(l-1)}=-u_0$  since otherwise 
$v_{\pi/(l-1)}-u_1$ would have  
$\nu(x_0)\ge l$ and thus more than $l$ nodal domains.
Since no intersection points outside $x_0$ are possible for functions 
in $U$, the nodal set $\Cal N(v_\al)$ moves smoothly to itself by 
this rotation. If the first nodal ray (counting from the angle 0) not 
leading to the outer boundary is 
connected by a smooth loop to the one with number $i$, 
then it follows that the second one is connected to the one with 
number $(i+1)$. But this is not possible without further intersection 
point, a contradiction. Thus $u_x$ is unique up to a multiplicative 
constant.

Finally we get a smooth mapping $D\to P(U)$, since $u_x$ is 
the solution of a linear system which has maximal rank by uniqueness. 
\qed\enddemo

\proclaim{\nmb.{3.8}. Lemma}
The mapping $x\mapsto u_x$ is continuous from $\overline D$ into the 
real projective space $P(U)$. 
\endproclaim

\demo{Proof}
Inside $D$ the function $u_x$, suitably normalized, depends smoothly 
on $x\in D$, by \nmb!{3.7}.
On the boundary $\p D$ the function $u_y$ depends smoothly on 
$y\in \p D$ by \nmb!{3.3}.

So it remains to show that $u_{x_n}\to u_y$ in the projective space $P(U)$ 
if the sequence $x_n$ in $D$ converges to $y\in \p D$. 
Since $P(U)$ is compact it suffices to show that each accumutation point 
of the sequence $u_{x_n}$ in $P(U)$ coincides with $u_y$. 

Thus let $v\in P(U)$ be a cluster point, then there is a subsequence 
of $u_{x_{n_k}}$ which converges to $v$. 

Let $C$ be a closed disk of small radius $\ep>0$ with center 
$y$ intersected with $\overline D$. 
Choose $n_k$ such that $x_{n_k}$ is still in 
the interior of $C$. Then of the $2l-2$ nodal rays of $u_{x_{n_k}}$ 
leaving $x_{n_k}$ all but one have to leave $C$, since otherwise 
there would exist a nodal domain which is completely contained in 
$C$. Since $\ep$ is small, this is not possible by energy reasons. 
\vskip 3mm 
\centerline{\epsfxsize=8cm\epsfbox{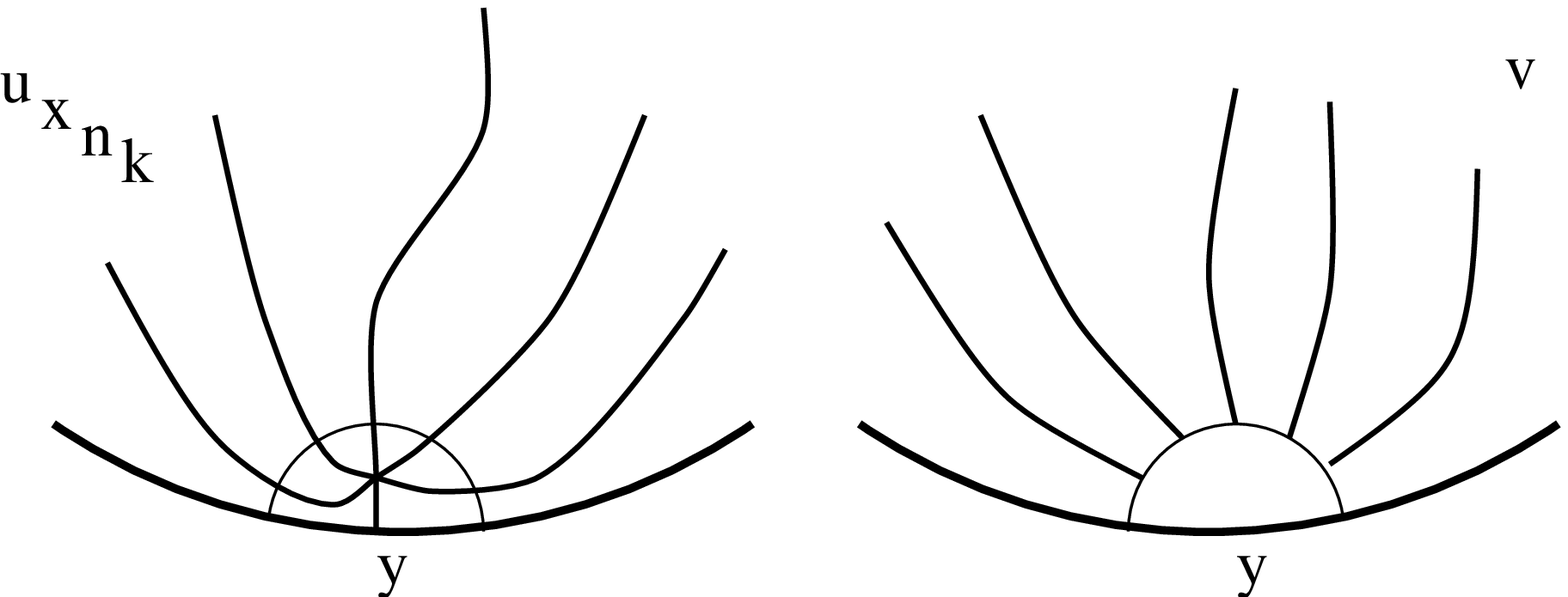}}
\medskip
But since $u_{x_{n_k}}$ converges to $v$ in $P(U)$, also at least 
$2l-3$ nodal lines of $v$ lead into $C$. Since $\ep$ was arbitrary, 
$2l-3$ nodal lines of $v$ hit $\p D$ at $y$. But the eigenfunction with 
this property is unique in $P(U)$ and is called $u_y$, by \nmb!{3.3}.
Thus $v=u_y$ in $P(U)$.
\qed\enddemo

\subhead{\nmb.{3.9} Proof of theorem B}\endsubhead
Suppose that $D$ is a simply connected domain and that $\p D$ is its
boundary. In \nmb!{3.8} we proved that the mapping $x\mapsto u_x$ is
continuous $\overline D\to P(U)$. 
Let $c:[0,2\pi]\to D$ be a closed smooth curve following $\p D$ 
anti-clockwise close enough so that all arguments below work. 
We want to analyze how the star of tangents at $c(t)$ to the
$l-1$ nodal lines of $u_{c(t)}$ crossing at $c(t)$ turns if we follow
$t$ from $0$ to $2\pi$. 

To make this precise, we consider the continuous function $f:D\to S^1$,
given by $f(x)=(2l-2)\al(x)$ modulo $2\pi$, where 
$$
t\mapsto t e^{2\pi i k/(2l-2) + i\al(x)},\quad k=0,\dots,2l-3, \quad t\ge 0
$$
are the tangents rays of the nodal lines through $x$ in $\Cal N(u_x)$. 
We want to analyze $f(c(t))$.

We consider again the sets the disjoint partition of $\p D$ into the
sets
$$\align
(\p D)_{2l-3} :&= \{y\in \p D: \rh_{u_y}(y)=2l-3\},\\
(\p D)_{2l-2} :&= \{y\in \p D: \rh_{u_y}(y)=2l-2\}.\\
\endalign$$
We note first the following fact:
\roster
\item If $x\in D$ is near enough to some point $y$ in the open set 
	$(\p D)_{2l-3}$ then the nodal pattern $\Cal N(u_x)$ can be read off
	the the nodal pattern $\Cal N(u_y)$, since the nodal domains move
	continuously. The nodal line leaving $\p D$ vertically stays
	connected to $\p D$ near $y$. For example:
\vskip 3mm 
\centerline{\epsfxsize=8cm\epsfbox{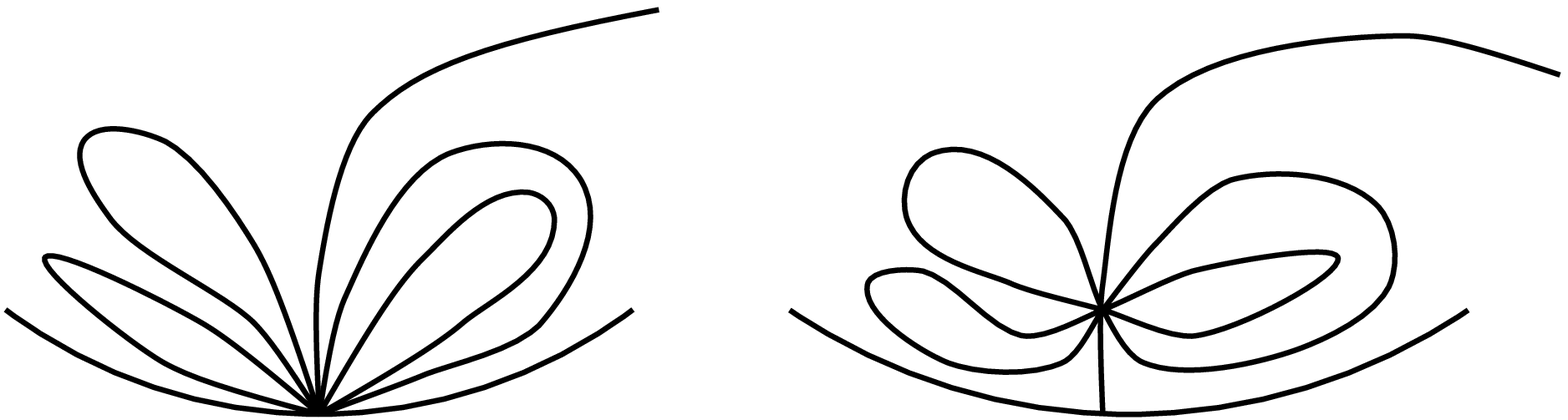}}
\medskip
\endroster
This is seen as follows: $\Cal N(u_x)$ for $x\in D$ can have at most
two nodal lines connecting $x$ to $\p D$, since otherwise there would
be too many nodal domains by \nmb!{2.8}. Moreover nodal lines can move
off $\p D$ only in pairs. 
Thus, if $x$ moves from $y\in (\p D)_{2l-3}$ into $D$, the nodal 
lines of $\Cal N(u_x)$ 
move away from $\p D$ in pairs and one of them stays connected to $\p D$,
since there was an odd number of them at $y$. Let us call this nodal 
arc from $x$ to $\p D$ the {\it short arc} of $\Cal N(u_x)$: it 
exists if $x$ is near $(\p D)_{2l-3}$.
Since loops have to stay loops, the result follows.

This already implies that $f(c(t))$ follows the angle of $\p D$ along
each arc in $(\p D)_{2l-3}$. 

What happens at a point in $(\p D)_{2l-2}$? Without loss we assume
that this point is $0\in (\p D)_{2l-2}$, that $\p D$ has horizontal
tangent at 0, and that $D$ lies above. Then there are two connected
components of $(\p D)_{2l-3}$ to the left and to the right of $0$: 
the open arcs $I_1$ and $I_2$. 
We claim that:
\roster
\item[2] When passing over $0$ from above $I_1$ to above $I_2$, the 
       angle $\al(c(t))$ increases by an amount of $2\pi/(2l-2)$. 
\endroster
To see this,
from \nmb!{3.6} we conclude that for
$y$ to the left, in $I_1$,
the last hitting point $z(y)$ of the nodal
pattern $\Cal N(u_y)$ has to move through $I_2$ towards $0$ if $y$
moves right to $0$; if $y$ then continues to move right, away from 0, 
then the leftmost nodal line emanating from $0$ has to move off to 
the left and to slide away through $I_1$. 
\vskip 3mm 
\centerline{\epsfxsize=9cm\epsfbox{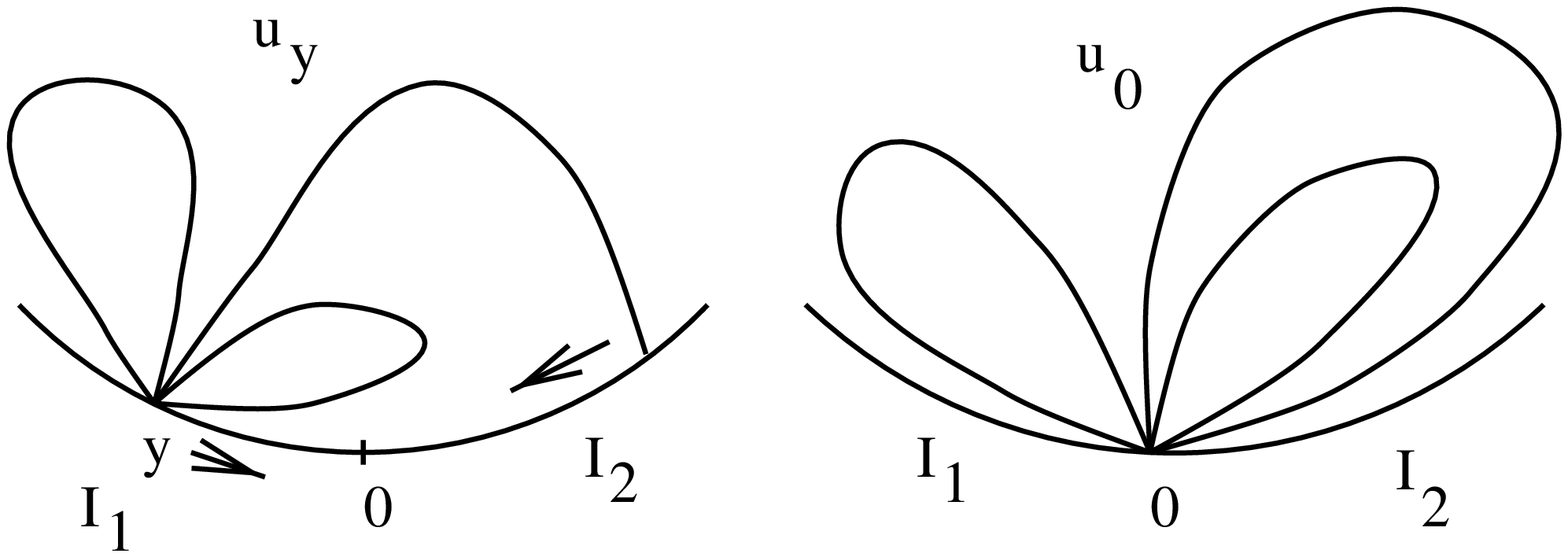}}
\medskip
If $c(t)$ moves from the left over 0 the nodal pattern 
$\Cal N(u_{c(t)})$ follows closely the behavious above. Well to the 
left, two nodal arcs connect $c(t)$ to $\p D$, the short one directly 
to $I_1$ below, and another one, call it the long one, far away, but 
moving through $I_2$ towards the short one. 

Eventually, near 0, they have to meet, to lift off $\p D$, and the 
next loop in clockwise direction of the short arc has to touch 
$\p D$, so that the short arc is replace by its next neighbor in 
clockwise direction and the new long arc then has to move away 
through $I_1$. Namely, this is the only continuous behaviour which 
connects the behaviour to the left to the one at the right of $0$, 
where we know exactly what happens. 
An illustration of the behaviour is the following:
\vskip 3mm 
\centerline{\epsfxsize=10cm\epsfbox{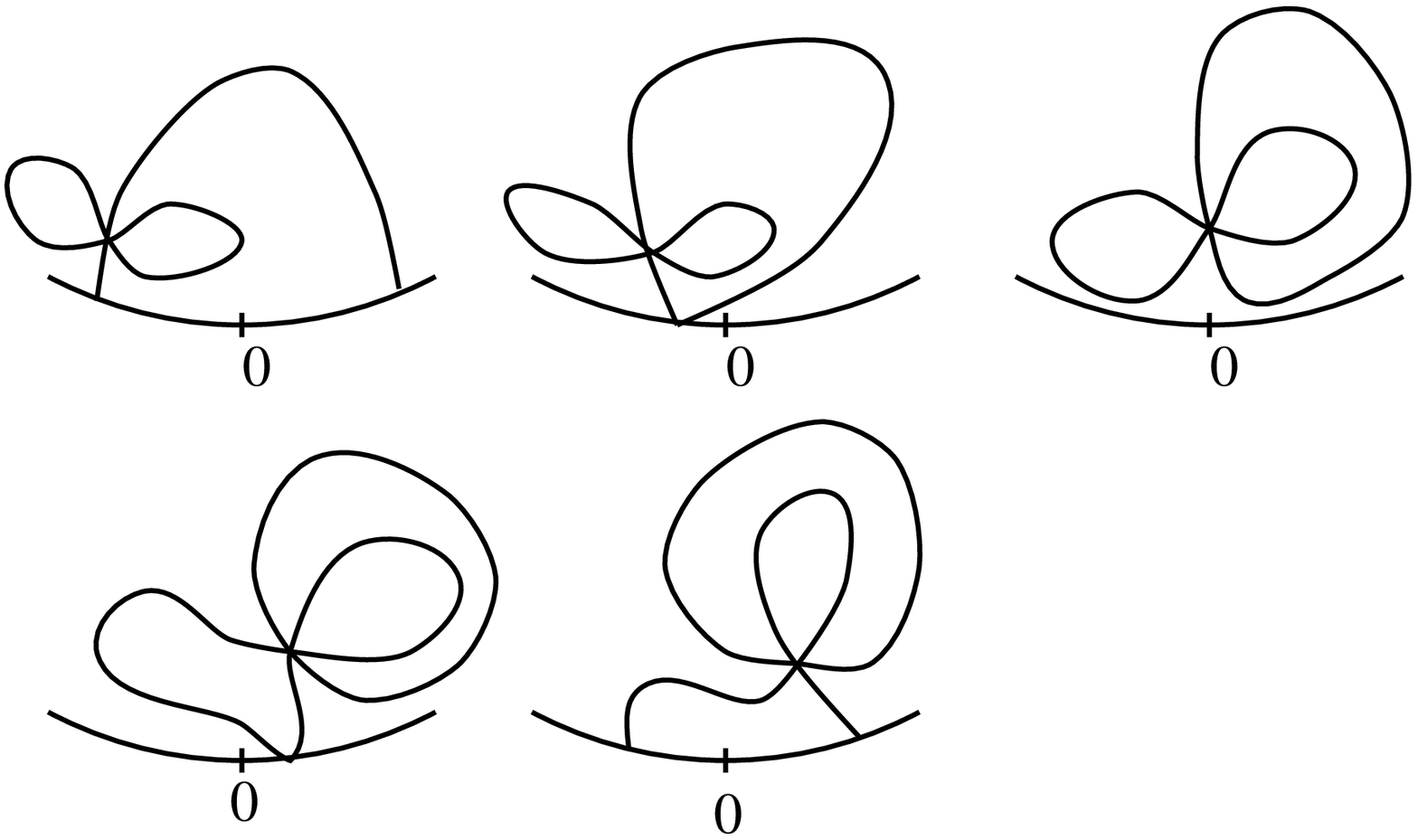}}
\medskip
Thus statement \therosteritem2 is proved. 

So finally the smooth mapping $t\mapsto f(c(t))$, $S^1\to S^1$ has
mapping degree $2\pi\#(\p D)_{2l-2}+ 2\pi(2l-2)>0$ and cannot be null
homotopic. But by construction it is continuously extended into the
interior of the circle and thus is nullhomotopic, a contradiction.
This finishes the proof for the simply connected case.

If $D$ is not simply connected, let $(\p D)^i$ for $i=1,\dots,p$ be the 
connected components of $\p D$. Choose a point $x_1$ near $(\p D)^1$. 
Then we choose a smooth curve $c:S^1\to D$ which starts at $x^1$ and 
follows $(\p D)^1$ closely back to $x_1$, then from $x_1$ along a 
smooth path $e_2$ to a point $x_2$ near $(\p D)^2$, then follows 
$(\p D)^2$ closely back to $x_2$, then back along $e_2$ to $x_1$. Then 
it follows a path $e_3$ not intersecting $e_2$ to some point $x_3$ 
and $\d D^3$, and so on until we end again at $x_1$. 
\vskip 3mm 
\centerline{\epsfxsize=10cm\epsfbox{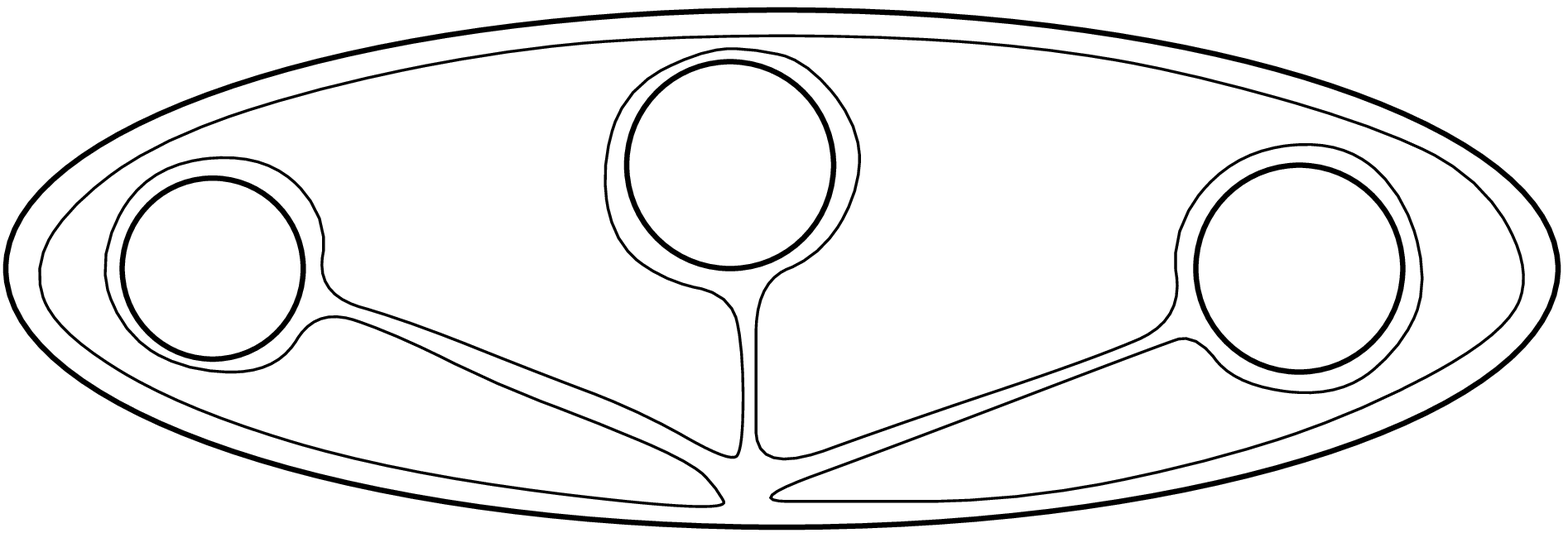}}
\medskip
Note that all results above also work for non simply connected 
domains, since we always worked with eigenfunctions 
which have the maximal number of nodal domains allowed by \nmb!{2.8}:
each of the further (inner) boundary components can be hit by at most 
one nodal line twice, otherwise we get too many 
nodal domains. 

Furthermore, all boundary components are equivalent for 
our arguments (put $D$ into $S^2$), and we treat each of them 
separately.

We consider again $f(c(t))$. 
At each boundary component the 
contribution to the mappingg degree of $f$ is a positive integer, by 
the arguments given above. The contributions from the parts going 
along the $e_i$ cancel each other. So the mapping $f\o c: S^1\to S^1$ 
has positive mapping degree, and thus cannot be nullhomotopic. 
But the curve $c(t)$ bounds a simply connected region, thus the 
mapping $f\o c:S^1\to S^1$ has a continuous extension $f$ to the 2 
cell in the interior. So it is nullhomotopic, a contradiction.
\qed

\enddocument